\documentclass[runningheads]{llncs}
\def\zibreport{1}

\usepackage[breaklinks, colorlinks, linkcolor=blue, citecolor=blue, urlcolor=black]{hyperref}

\usepackage[utf8]{inputenc}
\usepackage{graphicx}
\usepackage{todonotes}
\usepackage{tabularx}
\usepackage{tabularcalc}
\usepackage{booktabs}
\usepackage{multirow}
\usepackage{algorithm}
\usepackage{algpseudocode}
\usepackage{subfig}
\usepackage{siunitx}
\usepackage{zibtitlepage}
\usepackage{hyperref}
\usepackage{amsmath}
\usepackage{amssymb}
\usepackage{xspace}
\usepackage{xcolor}
\usepackage{collcell}

\newcommand{\eg}{e.g.,\xspace}

\newcommand{\bandb}{branch-and-bound\xspace}

\newcommand{\Abs}[1]{\vert #1\vert}

\newcolumntype{R}{>{\collectcell\ApplyColor}{r}<{\endcollectcell}}
\newcommand{\solver}[1]{\textsc{#1}\xspace}

\newcommand{\scip}{\solver{SCIP}}

\newcommand{\soplex}{\solver{SoPlex}}
\newcommand{\soplexversion}{5.0.2}
\newcommand{\soplexv}{\solver{SoPlex}~\soplexversion\xspace}
\newcommand{\papilo}{\solver{PaPILO}}

\newcommand{\qsoptex}{\solver{QSopt\_ex}}
\newcommand{\cplex}{\solver{CPLEX}}

\newcommand{\miplib}{\mbox{MIPLIB}}

\newcommand{\Z}{\mathbb{Z}\xspace}

\newcommand{\Q}{\mathbb{Q}\xspace}

\newcommand{\I}{\mathcal{I}\xspace}

\newcommand{\percentage}[2]{\pgfmathparse{#1/#2 *100}\pgfmathprintnumber[precision=1]{\pgfmathresult}\%}
\newcommand{\reduction}[2]{\pgfmathparse{(#2-#1)/#2 *100}\pgfmathprintnumber[precision=1]{\pgfmathresult}\%}
\newcommand{\increase}[2]{\pgfmathparse{(#1-#2)/#2 *100}\pgfmathprintnumber[precision=1]{\pgfmathresult}\%}
\newcommand{\fraction}[2]{\pgfmathparse{(#1)/#2}\pgfmathprintnumber[precision=1]{\pgfmathresult}}

\newcommand{\colorquot}[2]{
	\pgfmathparse{(#1<0.9*#2)?1:0}\ifdim\pgfmathresult pt>0pt \textcolor{blue}{$\mathbf{\reduction{#1}{#2}}$}\else
	\pgfmathparse{(#1>1.1*#2)?1:0}\ifdim\pgfmathresult pt>0pt \textcolor{red}{$\mathbf{\reduction{#1}{#2}}$} \else
	$\reduction{#1}{#2}$\fi \fi
}

\newcommand{\quotient}[2]{\pgfmathparse{(#1)/#2 }\pgfmathprintnumber[precision=1]{\pgfmathresult}}

\newcommand{\easy}{\textsc{fpeasy}\xspace}
\newcommand{\numdiff}{\textsc{numdiff}\xspace}

\begin{document}
\ifthenelse{\zibreport = 1}{

\newcommand{\tinytodo}[1]{\todo[size=\small]{#1}\xspace}
\newcommand{\myorcidlink}[1]{\,\href{https://orcid.org/#1}{\raisebox{-0.45ex}{\includegraphics[width=1.8ex]{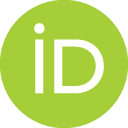}}}}

\setlength{\marginparwidth}{4.2cm}

\ZTPTitle{A Computational Status Update for Exact Rational Mixed Integer Programming}
\ZTPAuthor{
   Leon~Eifler\protect\myorcidlink{0000-0003-0245-9344},
   Ambros~Gleixner\protect\myorcidlink{0000-0003-0391-5903}
}

\ZTPNumber{21-04}
\ZTPMonth{January}
\ZTPYear{2021}

\ZTPInfo{This paper has been accepted to IPCO 2021.
Please cite as: \\ \textit{Leon Eifler, Ambros Gleixner, A computational Status Update for Exact Rational Mixed Integer Programming. 
\\ Accepted for publication in Integer Programming and Combinatorial Optimization: 22nd International Conference, IPCO 2021}}

\zibtitlepage

}{}


\title{A Computational Status Update for Exact Rational Mixed Integer Programming\thanks{The work for this article has been conducted within the Research Campus Modal funded by the German Federal Ministry of Education and Research (BMBF grant numbers 05M14ZAM, 05M20ZBM).}}
\author{Leon Eifler\inst{1}\orcidID{0000-0003-0245-9344} \and
  Ambros Gleixner\inst{1,2}\orcidID{0000-0003-0391-5903}
}
\titlerunning{A Computational Status Update for Exact Rational MIP}
\institute{%
  Zuse Institute Berlin, Takustr.~7, 14195 Berlin, Germany\\
  \email{\{eifler,gleixner\}@zib.de}
  \and
  HTW Berlin, Ostendstraße 1, 12459 Berlin, Germany
}

\maketitle


\begin{abstract}
The last milestone achievement for the roundoff-error-free solution of general
mixed integer programs over the rational numbers was a hybrid-precision
branch-and-bound algorithm published by Cook, Koch, Steffy, and Wolter in 2013.
We describe a substantial revision and extension of this framework that
integrates symbolic presolving, features an exact repair step for solutions from
primal heuristics, employs a faster rational LP solver based on LP iterative
refinement, and is able to produce independently verifiable certificates of
optimality.
We study the significantly improved performance and give insights into the
computational behavior of the new algorithmic components.
On the MIPLIB~2017 benchmark set, we observe an average speedup of
$\fraction{3317.5}{505.1}$x over the original framework and $\fraction{47}{17}$
times as many instances solved within a time limit of two hours.
\end{abstract}

\section{Introduction}

It is widely accepted that mixed integer programming (MIP) is a powerful tool for solving a broad variety of challenging optimization problems and that state-of-the-art MIP solvers are sophisticated and complex computer programs. However, virtually all established solvers today rely on fast floating-point arithmetic. Hence, their theoretical promise of global optimality is compromised by roundoff errors inherent in this incomplete number system. Though tiny for each single arithmetic operation, these errors can accumulate and result in incorrect claims of optimality for suboptimal integer assignments, or even incorrect claims of infeasibility.
Due to the nonconvexity of MIP, even performing an {a posteriori} analysis of
such errors or postprocessing them becomes difficult.


In several applications, these numerical caveats can become actual limitations.
This holds in particular when the solution of mixed integer programs is used as
a tool in mathematics itself.
Examples of recent work that employs MIP to investigate open mathematical
questions include
\cite{Bofi19,Burt12,eifler2020safe,KenterEtAl2018,LanciaEtAl2020,Pula20}.
Some of these approaches are forced to rely on floating-point solvers because
the availability, the flexibility, and most importantly the computational
performance of MIP solvers with numerically rigorous guarantees is currently
limited.
This makes the results of these research efforts not as strong as they could be.
Examples for industrial applications where the correctness of results is
paramount include hardware verification~\cite{Achterberg2007} or compiler
optimization~\cite{WilkenEtAl2000}.


The milestone paper by Cook, Koch, Steffy, and
Wolter~\cite{CookKochSteffyetal.2013} presents a hybrid-precision
branch-and-bound implementation that can still be considered the state of the
art for solving general mixed integer programs exactly over the rational
numbers.
It combines symbolic and numeric computation and applies different dual bounding
methods~\cite{Espinoza2006,Neumaier02safebounds,SteffyWolter2013} based on
linear programming (LP) in order to dynamically trade off their speed against
robustness and quality.



However, beyond advanced strategies for branching and bounding,
\cite{CookKochSteffyetal.2013} does not include any of the supplementary
techniques that are responsible for the strong performance of floating-point MIP
solvers today.
In this paper, we make a first step to address this research gap in two main
directions.

First, we incorporate a \emph{symbolic presolving} phase,
which safely reduces the size and tightens the formulation of the
instance to be passed to the branch-and-bound process.
This is motivated by the fact that presolving has been identified by several
authors as one of the components---if not \emph{the} component---with the largest
impact on the performance of floating-point MIP
solvers~\cite{AchterbergBixbyGuetal.2019,AchterbergWunderling2013}.
To the best of our knowledge, this is the first time that the impact of symbolic
preprocessing routines for general MIP is analyzed in the literature.

Second, we complement the existing dual bounding methods by enabling the use of
\emph{primal heuristics}.
The motivation for this choice is less to reduce the total solving time, but
rather to improve the usability of the exact MIP code in practical settings
where finding good solutions earlier may be more relevant than proving
optimality eventually.
Similar to the dual bounding methods, we follow a hybrid-precision scheme.
Primal heuristics are exclusively executed on the floating-point approximation
of the rational input data.
Whenever they produce a potentially improving solution, this solution is checked
for approximate feasibility in floating-point arithmetic.
If successful, the solution is postprocessed with an exact repair step that
involves an exact LP solve.

Moreover,
we integrate the exact LP solver SoPlex, which follows the recently developed
scheme of \emph{LP iterative refinement}~\cite{GleixnerSteffy2020},
we extend the logging of \emph{certificates} in the recently developed VIPR
format to all available dual bounding methods~\cite{VIPR},
and produce a thoroughly \emph{revised implementation} of the original
framework~\cite{CookKochSteffyetal.2013}, which improves multiple technical details.
Our computational study evaluates the performance of the new algorithmic aspects
in detail and indicates a significant overall speedup compared to the original
framework.

The overarching goal and contribution of this research is to extend the
computational practice of MIP to the level of rigor that has been achieved in
recent years, for example, by the field of satisfiability
solving~\cite{WetzlerHeuleHunt2014}, while at the same time retaining most of
the computational power embedded in floating-point solvers.
In MIP, a similar level of performance and rigor is certainly much more
difficult to reach in practice, due to the numerical operations that are
inherently involved in solving general mixed integer programs.
However, we believe that there is no reason why this vision should be
\emph{fundamentally} out of reach for the rich machinery of MIP techniques
developed over the last decades.
The goal of this paper is  to demonstrate the viability of this agenda
within a first, small selection of methods.
The resulting code is freely available for research purposes as an extension of
\scip~7.0~\cite{GithubExactSCIP}.


\section{Numerically exact mixed integer programming}
\label{sec:background}



In the following, we describe related work in numerically exact
optimization, including the main ideas and features of the framework that we
build upon.
%
%
Before turning to the most general case, we would like to mention that
roundoff-error-free methods are available for several specific classes of pure
integer problems.
One example for such a combinatorial optimization problem is the traveling
salesman problem, for which the branch-and-cut solver Concorde applies safe
interval-arithmetic to postprocess LP relaxation solutions and ensures the
validity of domain-specific cutting planes by their combinatorial
structure~\cite{applegate2006concorde}.

A broader class of such problems, on binary decision variables, is addressed in
\emph{satisfiability solving} (SAT) and \emph{pseudo-Boolean optimization}
(PBO)~\cite{BiereHeulevanMaarenWalsh2009}.
Solvers for these problem classes usually do not suffer from numerical errors
and often support solver-independent verification of
results~\cite{WetzlerHeuleHunt2014}.
While optimization variants exist, the development of these methods is to a
large extent driven by feasibility problems.
The broader class of solvers for \emph{satisfiability modulo theories} (SMT), e.g.,
\cite{deMouraEtAl2008}, may also include real-valued variables, in particular
for satisfiability modulo the theory of linear arithmetic.
However, as pointed out also
in~\cite{FaureNieuwenhuisOliverasRodriguez-Carbonell2008}, the target
applications of SMT solvers differ significantly from the motivating use cases
in LP and MIP.

Exact optimization over convex polytopes intersected with lattices is also
supported by some software libraries for polyhedral
analysis~\cite{Polymake2017,BagnaraHZ08SCP}.
These tools are not particularly targeted towards solving LPs or MIPs of larger
scale and usually follow the naive approach of simply executing all operations
symbolically, in exact rational arithmetic.
This yields numerically exact results and can even be highly efficient
as long as the size of problems or the encoding length of intermediate numbers
is limited.
However, as pointed out by~\cite{Espinoza2006}
and~\cite{CookKochSteffyetal.2013}, this \emph{purely symbolic approach} quickly
becomes prohibitively slow in general.

By contrast, the most effective methods in the literature rely on a \emph{hybrid
  approach} and combine exact and numeric computation.
For solving pure LPs exactly, the most recent methods that follow this paradigm
are \emph{incremental precision boosting}~\cite{APPLEGATE2007} and \emph{LP iterative
refinement}~\cite{GleixnerSteffy2020}.
In an exact MIP solver, however, it is not always necessary to solve LP
relaxations completely, but it often suffices to provide dual bounds that
underestimate the optimal relaxation value safely.
This can be achieved by postprocessing approximate LP solutions.
\emph{Bound-shift}~\cite{Neumaier02safebounds} is such a method that only
relies on directed rounding and interval arithmetic and is therefore very fast.
However, as the name suggests it requires upper and lower bounds on all
variables in order to be applicable.
A more widely applicable bounding method is
\emph{project-and-shift}~\cite{SteffyWolter2013}, which uses an interior point
or ray of the dual LP.
These need to be computed by solving an auxiliary LP exactly in advance, though
only once per MIP solve.
Subsequently, approximate dual LP solutions can be corrected by projecting them
to the feasible region defined by the dual constraints and shifting the result
to satisfy sign constraints on the dual multipliers.

The hybrid branch-and-bound method of~\cite{CookKochSteffyetal.2013} combines such safe dual bounding methods with a state-of-the-art branching
heuristic, reliability branching~\cite{ACHTERBERG200542}.
It maintains both the exact problem formulation
\begin{align*}
   \min \{&c^T x \ \vert \ Ax \ge b,\ x \in \Q^n, \ x_i \in \Z \ \forall i \in \I\}
\end{align*}
with rational input data $A\in \Q^{m\times n}, c \in Q^n, b \in \Q^m$, as well
as a floating-point approximation with data $\bar A, \bar b, \bar c$, which are
defined as the componentwise closest numbers representable in floating-point
arithmetic.
The set $\I \subseteq \{1,\ldots,n\}$ contains the indices of integer variables.

During the solve, for all LP relaxations, the floating-point approximation is first solved in floating-point arithmetic as an approximation and then postprocessed to generate a valid dual bound.
The methods available for this safe bounding step are the previously described
\emph{bound-shift} \cite{Neumaier02safebounds}, \emph{project-and-shift}
\cite{SteffyWolter2013}, and an \emph{exact LP solve} with the exact LP solver
\qsoptex based on incremental precision boosting~\cite{APPLEGATE2007}.
(Further dual bounding methods were tested, but reported as less important
in~\cite{CookKochSteffyetal.2013}.)
On the primal side, all solutions are checked for feasibility in exact arithmetic before being accepted.

Finally, this exact MIP framework was recently extended by the possibility to
generate a \emph{certificate} of correctness~\cite{VIPR}. This certificate is a
tree-less encoding of the \bandb search, with a set of dual multipliers to prove
the dual bound at each node or its infeasibility.
Its correctness can be verified independently of the solving process using the
checker software VIPR~\cite{VIPRweb}.

\section{Extending and improving an exact MIP framework}
\label{sec:advances}

The exact MIP solver presented here extends \cite{CookKochSteffyetal.2013} in
four ways:
the addition of a symbolic presolving phase,
the execution of primal floating-point heuristics coupled with an exact repair
step,
the use of a recently developed exact LP solver based on LP iterative
refinement,
and a generally improved integration of the exact solving routines into the core
branch-and-bound algorithm.
\medskip

\noindent
\textbf{Symbolic presolving.}
The first major extension is the addition of symbolic presolving.
%
To this end, we integrate the newly available presolving library \papilo \cite{Papilo} for integer and linear programming.
\papilo has several benefits for our purposes.

First, its code base is by design fully templatized with respect to the arithmetic type.
This enables us to integrate it with rational numbers as data type for storing the MIP data and all its computations.
Second, it provides a large range of presolving techniques already implemented.
The ones used in our exact framework are
coefficient strengthening,
constraint propagation,
implicit integer detection,
singleton column detection,
substitution of variables,
simplification of inequalities,
parallel row detection,
sparsification,
probing,
dual fixing,
dual inference,
singleton stuffing,
and dominated column detection.
For a detailed explanation of these methods, we refer to \cite{AchterbergBixbyGuetal.2019}.
Third, \papilo comes with a sophisticated parallelization scheme that helps to
compensate for the increased overhead introduced by the use of rational arithmetic.
For details see~\cite{GamrathAndersonBestuzhevaetal.2020}.

When \scip enters the presolving stage, we pass a rational copy of the problem to
\papilo, which executes its presolving routines iteratively until no
sufficiently large reductions are found.
Subsequently, we extract the postsolving information provided by \papilo to transfer the model reductions to \scip.
These include fixings, aggregations, and bound changes of variables and strengthening or deletion of constraints, all of which are performed in rational arithmetic.
\medskip

\noindent
\textbf{Primal heuristics.}
The second extension is the safe activation of \scip's floating-point heuristics and the addition of an exact repair heuristic for their approximate solutions.
Heuristics are not known to reduce the overall solving time drastically, but they can be particularly useful on hard instances that cannot be solved at all, and in order to avoid terminating without a feasible solution.

In general, activating \scip's floating-point heuristics does not interfere with
the exactness of the solving process, although care has to be taken that no
changes to the model are performed, \eg the creation of a no-good constraint.
However, the chance that these heuristics find a solution that is feasible in the exact sense can be low, especially if equality constraints are present in the model.
Thus, we postprocess solutions found by floating-point heuristics in the following way. First, we fix all integer variables to the values found by the floating-point heuristic, rounding slightly fractional values to their nearest integer. Then an exact LP is solved for the remaining continuous subproblem. If that LP is feasible, this produces an exactly feasible solution to the mixed integer program.

Certainly, frequently solving this subproblem exactly can create a significant overhead compared to executing a floating-point heuristic alone, especially when a large percentage of the variables is continuous and thus cannot be fixed. Therefore, we impose working limits on the frequency of running the exact repair heuristic, which are explained in more detail in Sec.~\ref{sec:computations}.
\medskip

\noindent
\textbf{LP iterative refinement.}
Exact linear programming is a crucial part of the exact MIP solving process.
Instead of \qsoptex, we use \soplex as the exact linear programming solver. The reason for this change is that \soplex uses LP iterative refinement~\cite{GleixnerSteffyWolter2016} as the strategy to solve LPs exactly, which compares favorably against incremental precision boosting~\cite{GleixnerSteffy2020}.


\medskip

\noindent
\textbf{Further enhancements.}
We improved several algorithmic details in the implementation of the hybrid
branch-and-bound method. We would like to highlight two examples for these changes.
First, we enable the use of an \emph{objective limit} in the floating-point LP solver, which was not possible in the original framework.
Passing the primal bound as an objective limit to the floating-point LP solver
allows the LP solver to stop early just after its dual bound exceeds the global
primal bound.
However, if the overlap is too small, postprocessing this LP
solution with safe bounding methods can easily lead to a dual bound that no
longer exceeds the objective limit.
For this reason, before installing the primal bound as an objective limit in the
LP solver, we increase it by a small amount computed from the statistically
observed bounding error so far.
Only when safe dual bounding fails, the objective limit is solved again without
objective limit.

Second, we reduce the time needed for checking exact feasibility of primal
solutions by prepending a safe floating-point check.
Although checking a single solution for feasibility is fast, this happens often throughout the solve and doing so repeatedly in exact arithmetic can become computationally expensive.
To implement such a safe floating-point check, we employ \emph{running error analysis}~\cite{Higham2002}.
Let $x^* \in \Q^n$ be a potential solution 
and let $\bar x^*$ be the floating-point approximation of $x^*$.
Let $a \in \Q^n$ be a row of $A$ with floating-point approximation $\bar a$, and right hand side $b_j \in \Q$.
Instead of computing $\sum_{i=1}^n a_ix_i^*$ symbolically, we instead compute $\sum_{i=1}^n \bar a_i \bar x_i^*$ in floating-point arithmetic, and alongside compute a bound on the maximal rounding error that may occur. We adjust the running error analysis described in~\cite[Alg.~3.2]{Higham2002} to also account for roundoff errors $\Abs{\bar x^*-x^*}$ and $\Abs{\bar a-a}$.
After doing this computation, we can check if either $s - \mu \ge b_j$ or $s + \mu \le b_j$. In the former, the solution $x^*$ is guaranteed to fulfill $\sum_{i=1}^n a_i x^*_i \ge b_j$; in the latter, we can safely determine that the inequality is violated; only if neither case occurs, we recompute the activity in exact arithmetic.

We note that this could alternatively be achieved by directed rounding, which
would give tighter error bounds at a slightly increased computational effort.
However, empirically we have observed that most equality or inequality
constraints are either satisfied at equality, where an exact arithmetic check
cannot be avoided, or they are violated or satisfied by a slack larger
than the error bound $\mu$, hence the running error analysis is sufficient to
determine feasibility.



\section{Computational Study}
\label{sec:computations}


We conduct a computational analysis to answer three main questions.
\emph{First, how does the revised branch-and-bound framework compare to the previous implementation, and to which components can the changes be attributed?}
To answer this question, we compare the original framework \cite{CookKochSteffyetal.2013} against our improved implementation, including the exact LP solver \soplex, but with primal heuristics and exact presolving still disabled. In particular, we analyze the importance and performance of the different dual bounding methods.

\emph{Second, what is the impact of the new algorithmic components symbolic presolving and primal heuristics?}
To answer this question, we compare their impact on the solving time and the number of solved instances, as well as present more in-depth statistics, such as \eg the primal integral \cite{Berthold2013} for heuristics or the number of fixings for presolving.
In addition, we compare the effectiveness of performing presolving in rational and in floating-point arithmetic.

\emph{Finally, what is the overhead for producing and verifying certificates?} Here, we consider running times for both the solver and the certificate checker, as well as the overhead in the safe dual bounding methods introduced through enabling certificates.
This provides an update for the analysis in \cite{VIPR}, which was limited to the two bounding methods project-and-shift and exact LP.

The experiments were performed on a cluster of
Intel Xeon CPUs E5-2660 with 2.6\,GHz and 128\,GB main memory. 
As in \cite{CookKochSteffyetal.2013}, we use
\cplex as floating-point LP solver. {Due to compatibility issues, we needed to use \cplex 12.3.0 for the original and \cplex 12.8.0 for the new framework.
Although these versions are different, they are only used to solve floating-point LPs and have limited impact on the reported results: The vast majority of solving time is spent in the safe dual bounding methods.}
For exact LP solving, we use the same \qsoptex version as in \cite{CookKochSteffyetal.2013} and \soplex 5.0.2.
For all symbolic computations, we use the GNU Multiple Precision Library (GMP) 6.1.4 \cite{GMP}.
For symbolic presolving, we use \papilo 1.0.1 \cite{GamrathAndersonBestuzhevaetal.2020,Papilo}; all other SCIP presolvers are disabled.

As main test sets, we use the two test sets specifically curated in~\cite{CookKochSteffyetal.2013}:
one set with $57$~instances that were found to be easy for an inexact floating-point branch-and-bound solver (\easy), and one set of $50$~instances 
that were found to be numerically challenging, \eg due to poor conditioning or large coefficient ranges (\numdiff).
For a detailed description of the selection criteria, we refer to \cite{CookKochSteffyetal.2013}.
To complement these test sets with a set of more ambitious and recent instances, we conduct a final comparison on the \miplib~2017~\cite{Miplib2017} benchmark set.
%
All experiments to evaluate the new code are run with three different random
seeds, where we treat each instance-seed combination as a single observation.
As this feature is not available in the original framework, all comparisons with the original framework were performed with one seed. 
The time limit was set to 7200~seconds for all experiments.
If not stated otherwise all aggregated numbers are shifted geometric means with a shift of $0.001$~seconds or $100$ branch-and-bound nodes, respectively.

\bigskip

\noindent
\textbf{The branch-and-bound framework.}
As a first step, we compare the behavior of the safe branch-and-bound implementation from \cite{CookKochSteffyetal.2013}
with \qsoptex as the exact LP solver, against its revised implementation
with \soplexv as exact LP solver.
The original framework uses the ``Auto-Ileaved'' bounding strategy as recommended in \cite{CookKochSteffyetal.2013}. It dynamically chooses the dual bounding method, attempting to employ bound-shift as often as possible. An exact LP is solved whenever a node would be cut off within tolerances, but not with the exact the safe dual bound computed. In the new implementation we use a similar strategy, however we solve the rational LP relaxation every $5$~depth levels of the tree, due to improved performance in the exact LP solver.

\begin{figure}[t]
   \centering
   \hfill
   \begin{minipage}[t]{0.48\textwidth}
     \includegraphics[width=.98\textwidth]{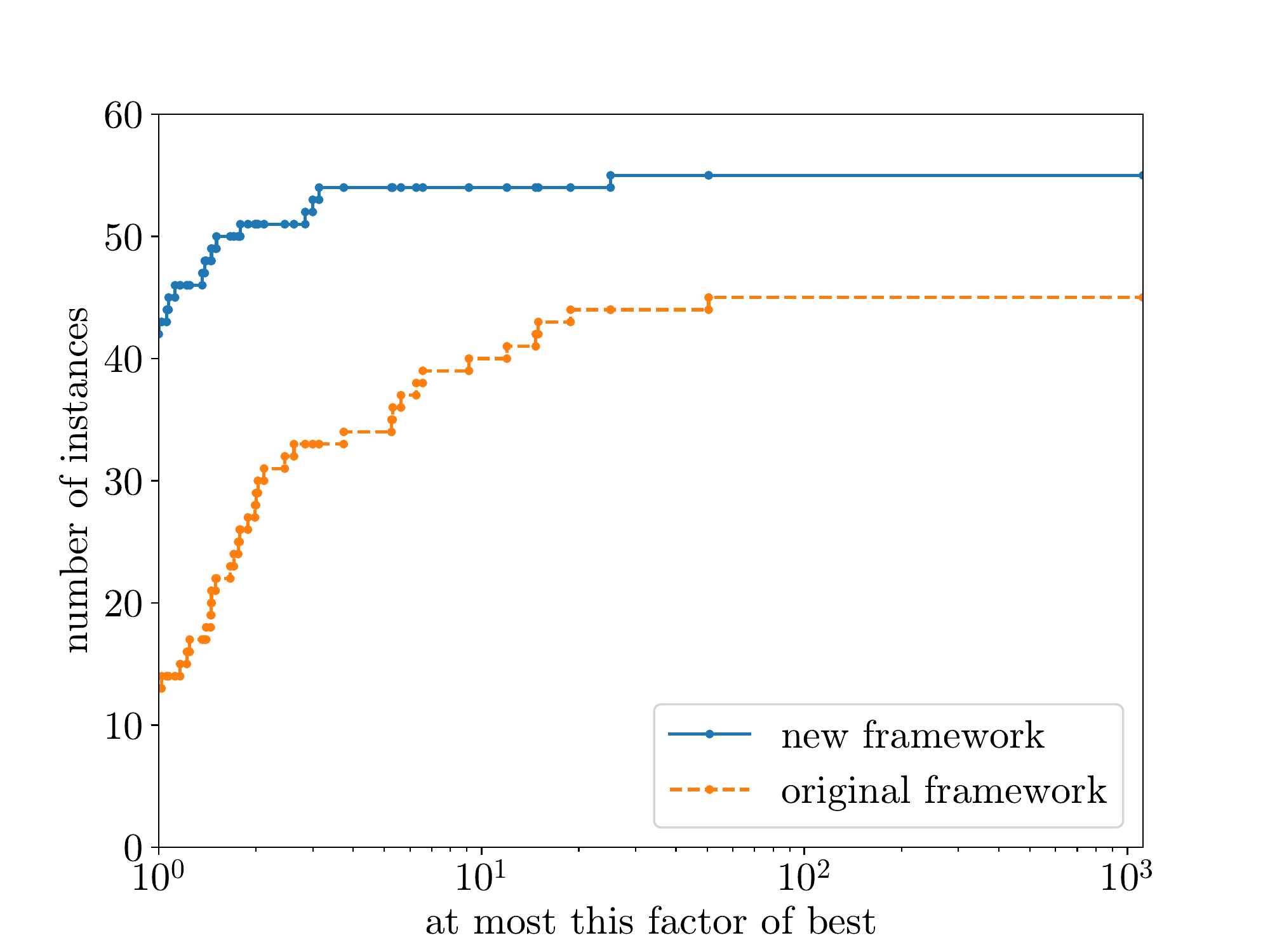}
  \end{minipage}
  \hfill
  \begin{minipage}[t]{0.48\textwidth}
     \includegraphics[width=.98\textwidth]{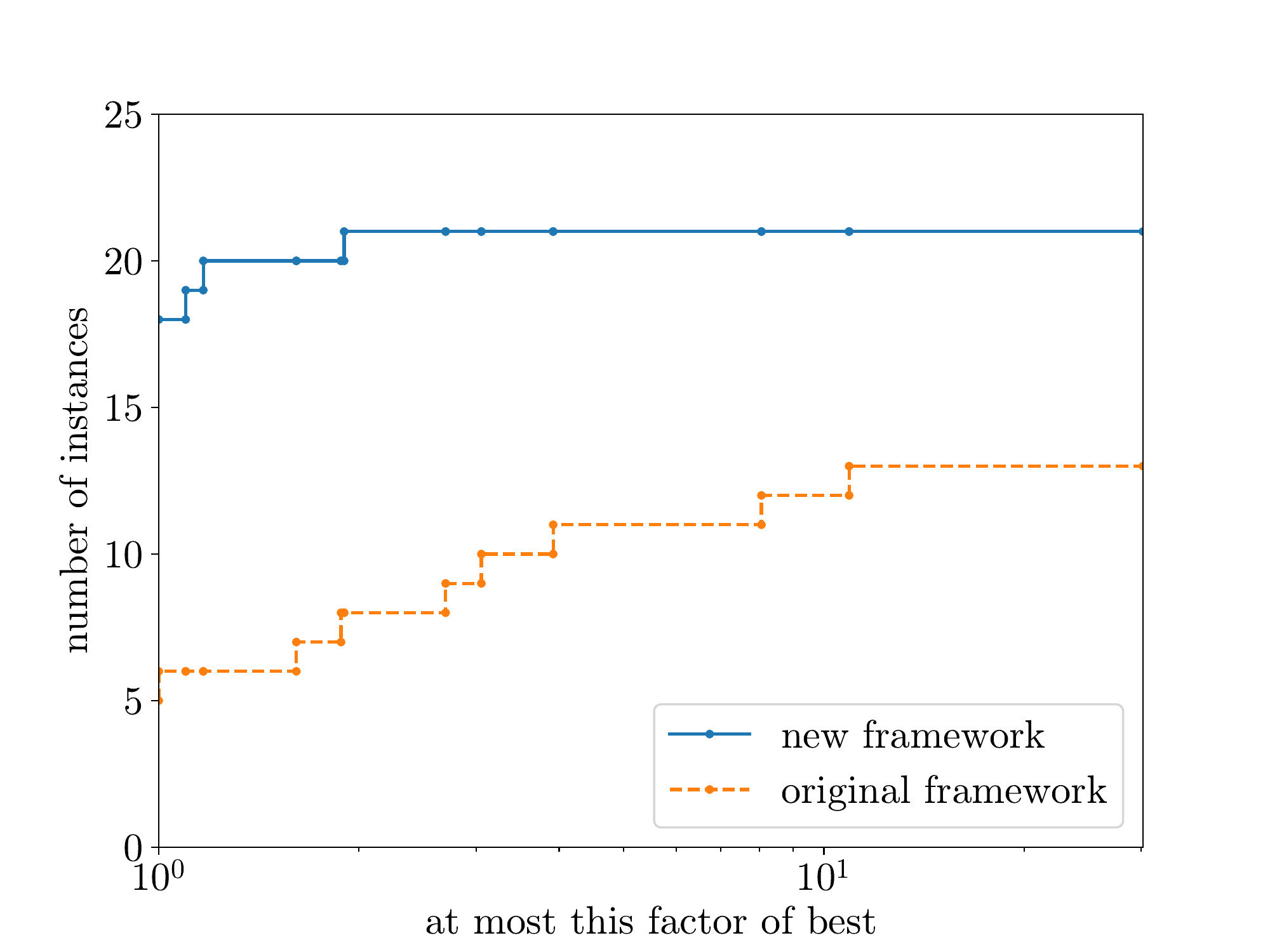}
  \end{minipage}
  \caption{Performance profiles comparing solving time of original and
       new framework without presolving and heuristics for \easy (left) and \numdiff (right)}
  \label{fig:pprof}
\end{figure}

Table \ref{oldvsnew-avg} reports the results for solving time, number of nodes, and total time spent in safe dual bounding (``dbtime''), for all instances that could be solved by at least one solver.
The new framework could solve $10$ instances more on \easy and $7$ more on \numdiff. On \easy, we observe a reduction of $\reduction{38.8}{128.4}$ in solving time and of $\reduction{11}{86.8}$ in safe dual bounding time.
On \numdiff, we observe a reduction of $\reduction{46.6}{237}$ in solving time, and of $\reduction{13.4}{114.6}$ in the time spent overall in the safe dual bounding methods.
We also see this significant performance improvement reflected in the two performance profiles in Fig.~\ref{fig:pprof}.

\begin{table}[t]
  \small{
    \caption{Comparison of original and new framework with presolving and primal heuristics disabled} 
    \medskip

     \label{oldvsnew-avg}
     \sisetup{round-mode=places,round-precision=1}
     \begin{tabular*}{\textwidth}{l@{\extracolsep{\fill}}r@{\extracolsep{\fill}}rS[table-format=4.1]S[table-format=5.1]S[table-format=4.1]rS[table-format=3.1]S[table-format=5.1]S[table-format=4.1]}
     \toprule
     & & \multicolumn{4}{c}{original framework} & \multicolumn{4}{c}{new framework} \\
      \cmidrule(){3-6}\cmidrule(){7-10}
      {Test set} & {size} & {solved} & \multicolumn{1}{r}{time} & \multicolumn{1}{r}{nodes} & {dbtime} & {solved} &  \multicolumn{1}{r}{time} & \multicolumn{1}{r}{nodes} & \multicolumn{1}{r}{dbtime} \\
     \midrule
     \easy & 55 & 45 & 128.42 &             8920.06 &        86.81 & 55 & 38.75 &             5940.79 &        10.97 \\ 
     \numdiff & 21 & 13 & 237.02 &             8882.67 &       114.62 & 20 & 46.59 &             6219.29 &        13.36 \\ 
     \midrule
   \end{tabular*}
     \vspace*{-1.2em}
  }
 \end{table}

\begin{table}[t]
   \caption{Comparison of safe dual bounding techniques}
   \medskip
   
     \label{oldvsnew-bounding}
     \small{
       \sisetup{round-mode=off}
       \begin{tabular*}{\textwidth}{ll@{\extracolsep{\fill}}S[table-format=1.4]S[table-format=1.4]S[table-format=1.4]S[table-format=1.4]S[table-format=1.4]S[table-format=1.4]}
     \toprule
      &  & \multicolumn{3}{c}{original framework} & \multicolumn{3}{c}{new framework} \\
      \cmidrule(){3-5}\cmidrule(){6-8}
     Test set & stats & \multicolumn{1}{l}{bshift\;\;} & {pshift\;\,} & \multicolumn{1}{r}{exlp} & {bshift\;\,} & \multicolumn{1}{r}{pshift} & \multicolumn{1}{r}{exlp}  \\
     \midrule
     {\easy} 
     & calls/node    & 0.92 &       0.44 &         0.28 &  0.53 &       0.39 &         0.06 \\ 
     & time/call [s]    & 0.0026 & 0.0022 & 0.050 & 0.0072 & 0.0066 & 0.010 \\ 
     & time/solving time  & {2.9\,\%} &    {40.3\,\%} &      {32.1\,\%} & {10.8\,\%} &    {27.8\,\%} &      {4.4\,\%} \\ 
     \midrule
     {\numdiff\;\;}
   & calls/node  & 0.78 &       0.36 &         0.52 & 0.36 &       0.39 &         0.28 \\ 
     & time/call [s] & 0.0055 & 0.0036 & 0.4197 & 0.0247 & 0.0356 & 0.1556 \\
     & time/solving time &  {1.4\,\%} &    {22.7\,\%} &      {62.2\,\%}  & {5.2\,\%} &    {24.8\,\%} &      {40.1\,\%}  \\
     \midrule
   \end{tabular*}
     \vspace*{-1.2em}
     }
 \end{table}

We identify a more aggressive use of project-and-shift and faster exact LP solves as the two key factors for this improvement.
In the original framework, project-and-shift is restricted to instances that had less than $10000$ nonzeros. One reason for this limit is that a large auxiliary LP has to be solved by the exact LP solver to compute the relative interior point in project-and-shift. With the improvements in exact LP performance, 
it proved beneficial to remove this working limit in the new framework.

The effect of this change can also be seen in the detailed analysis of bounding times
given in Table \ref{oldvsnew-bounding}. For calls per node and the fraction of bounding time per total solving time, which are normalized well, we report the arithmetic means; for time per call, we report geometric means over all instances where the respective bounding method was called at least once.

The fact that time per call for project-and-shift (``pshift'') in the new framework increased by a factor of $\quotient{66}{22}$ (\easy) and $\quotient{356}{36}$ (\numdiff) is for the reason discussed above---it is now also called on larger instances.
This is beneficial overall since it replaces many slower exact LP calls.
The decrease in exact LP solving time per call (``exlp'') by a factor of $\quotient{4197}{1556}$ (\numdiff) and $\quotient{50}{10}$ (\easy) can also partly be explained by this change, and partly by an overall performance improvement in exact LP solving due to the use of LP iterative refinement~\cite{GleixnerSteffyWolter2016}.
The increase in bound-shift time (``bshift'') is due to implementation details, that will be addressed in future versions of the code, but its fraction of the total solving time is still relatively low.
Finally, we observe a decrease in the total number of safe bounding calls per node.
One reason is that we now disable bound-shift dynamically if its success rate drops below $20\%$.

Overall, we see a notable speedup and more solved instances, mainly due to the better management of dual bounding methods and faster exact LP solving.


\bigskip

\noindent
\textbf{Symbolic presolving.}
Before measuring the overall performance impact of exact presolving, we address
the question how effective and how expensive presolving in rational arithmetic
is compared to standard floating-point presolving.
For both variants, we configured \papilo to use the same tolerances for
determining whether a reduction found is strong enough to be accepted.
The only difference in the rational version is that all computations are done in
exact arithmetic and the tolerance to compare numbers and the
feasibility tolerance are zero.
Note that a priori it is unclear whether rational presolving yields more or less reductions.
Dominated column detection may be less successful due to
the stricter comparison of coefficients; the dual inference presolver might be
more successful if it detects earlier that a dual multiplier is strictly bounded
away from zero.

Table~\ref{presol-fpvsex} presents aggregated results for presolving time, the
number of presolving rounds, and the number of found fixings, aggregations, and
bound changes.
We use a shift of $1$ for the geometric means of rounds, aggregations, fixings, and bound changes to account for instances where presolving found no such reductions.
Remarkably, both variants yield virtually the same results on \easy.
On \numdiff, there are small differences, with a slight decrease in the number of fixings and aggregations and a slight increase in the number of bound changes for exact variant. The time spent for exact presolving increases by more than an order of magnitude
but symbolic presolving is still not a performance bottleneck.
It consumed only $0.86\%$ (\easy) and $\percentage{1.23}{58.7}$ (\numdiff) of the total solving time, as seen in Table \ref{presol-avg}. 
Exploiting parallelism in presolving provided no measureable benefit for floating-point presolving, but reduced symbolic presolving time by $\reduction{0.14}{0.25}$ (\easy) to $\reduction{0.5}{0.89}$ (\numdiff).
However, this benefit can be more pronounced on individual instances, \eg on \texttt{nw04}, where parallelization reduces the time for rational presolving by a factor of $\fraction{1770}{277}$ from 1770 to 277~seconds.

To evaluate the impact of exact presolving, we compare the performance of the basic branch-and-bound algorithm established above against the performance with presolving enabled. The results for all instances that could be solved to optimality by at least one setting are presented in Table \ref{presol-avg}.
Enabling presolving solves $3$ more instances on \easy and $10$ more instances on \numdiff. We observe a reduction in solving time of $\reduction{25.5}{42.1}$ (\easy) 
and $\reduction{58.7}{216.6}$ (\numdiff). 
The stronger impact on \numdiff is correlated with the larger number of
reductions observed in Table~\ref{presol-fpvsex}.

\begin{table}[t]
  \caption{Comparison of exact and floating-point presolving}
  \medskip

  \label{presol-fpvsex}
   \small{
   \sisetup{round-mode=places,round-precision=1}
   \begin{tabular*}{\textwidth}{l@{\extracolsep{\fill}}lS[table-format=1.2,round-precision=2,table-number-alignment=right]S[table-format=2.1,table-number-alignment=right]S[table-format=2.1,table-number-alignment=right]S[table-format=2.1,table-number-alignment=right]S[table-format=2.1,table-number-alignment=right]S[table-format=1.2,round-precision=2,table-number-alignment=right]S[table-format=2.1,table-number-alignment=right]S[table-format=2.1,table-number-alignment=right]S[table-format=2.1,table-number-alignment=right]S[table-format=2.1,table-number-alignment=right]}
   \toprule
   & & \multicolumn{5}{c}{floating-point presolving} & \multicolumn{5}{c}{exact presolving} \\
   \cmidrule(){3-7}\cmidrule(){8-12}
   {Test set} & {thrds} & {time} & \multicolumn{1}{r}{rnds} & \multicolumn{1}{r}{fixed} & \multicolumn{1}{r}{agg} & {bdchg} & {time} & {rnds} & \multicolumn{1}{r}{fixed} & \multicolumn{1}{r}{agg} & {bdchg} \\
   \midrule
   \easy & 1 &  0.0107 &                3.15 &                 8.53 &              3.50 &               10.39 &  0.2500 &                3.16 &                 8.53 &              3.50 &               10.38 \\
    & 20 &  0.0102 &                3.15 &                 8.53 &              3.50 &               10.39 &  0.1437 &                3.16 &                 8.53 &              3.50 &               10.38 \\
   \midrule
   \numdiff & 1 & 0.0406  & 8.27 &                53.76 &             55.68 &               51.36 &  0.8893 & 7.17 &               41.39 &             42.86 &               55.84 \\ 
    & 20 & 0.0395 & 8.27 &                53.76 &             55.68 &               51.36 &    0.4997 & 7.17 &            41.39 &             42.86 &               55.84 \\ 
   \midrule
   \end{tabular*}
     \vspace*{-1.2em}
   }
 \end{table}

\begin{table}[t]
  \caption{Comparison of new framework with and without presolving (3 seeds)}
  \medskip

   \label{presol-avg}
   \small{
     \sisetup{round-mode=places,round-precision=1}
   \begin{tabular*}{\textwidth}{l@{\extracolsep{\fill}}r@{\extracolsep{\fill}}rS[table-format=3.1,table-number-alignment=right]S[table-format=4.1,table-number-alignment=right]rS[table-format=3.1,table-number-alignment=right]lS[table-format=4.1,table-number-alignment=right]}
   \toprule
   & & \multicolumn{3}{c}{presolving disabled} & \multicolumn{4}{c}{presolving enabled} \\
   \cmidrule(){3-5}\cmidrule(){6-9} 
   {Test set} & {size} & {solved} & {time} &  {nodes} &  {solved} & {time} & {\hspace*{-3ex}(presolving)} & {nodes} \\
   \midrule
   \easy & 168 & 165 & 42.14 &             6145.34 & 168 & 25.52 & \hspace*{-3ex}(0.22) &          4724.05 \\ 
   \numdiff & 91 & 66 &216.62 &             7237.17 & 86 &  58.68&             \hspace*{-3ex}(1.23) &          2867.19 \\ 
   \midrule
   \end{tabular*}
     \vspace*{-1.2em}
   }
\end{table}

\bigskip

\noindent
\textbf{Primal heuristics.}
To improve primal performance, we enabled all SCIP
heuristics that the floating-point version executes by default.
To limit the fraction of solving time for the repair heuristic described in Section \ref{sec:advances}, the repair heuristic is only allowed to run at any point in the solve, if it was called at most half as often as the exact LP calls for safe dual bounding.
Furthermore, the repair heuristic is disabled on instances with more than $80\%$ continuous variables, since the overhead of the exact LP solves can drastically worsen the performance on those instances.
Whenever the repair step is not executed, the floating-point solutions are checked directly for exact feasibility.

First, we evaluate the cost and success of the exact repair heuristic over all instances where it was called at least once. The results are presented in Table \ref{heur-analysis}. The repair heuristic is effective at finding feasible solutions with a success rate of $46.9\%$ (\easy) and $25.6\%$ (\numdiff).
%
The fraction of the solving time spent in the repair heuristic is well below $1\%$. Nevertheless, the strict working limits we imposed are necessary
since there exist outliers for which the repair heuristic takes more than $5\%$ of the total solving time, and performance on these instances would quickly deteriorate if the working limits were relaxed.

\begin{table}[t]
  \caption{Statistics of repair heuristic for instances where repair step was called}
  \medskip

  \label{heur-analysis}
   \small{
   \begin{tabular*}{\textwidth}{l@{\extracolsep{\fill}}r@{\extracolsep{\fill}}S[table-format=3.1,round-mode=places,round-precision=1,table-number-alignment=right]rrrr}
   \toprule
   & & \multicolumn{4}{c}{time} & \\
   \cmidrule(){3-6}
   {Test set} & {size} & {total solving} & {repair} & {fail} & {success} & {success rate} \\
   \midrule
   \easy & 82 &  39.76 &     0.0020 &                 0.0017 &                 0.0003  &            46.9\,\%  \\
   \numdiff & 42 &   383.63 &   0.0187 &                 0.0077 &                 0.0062 &               25.6\,\%  \\
   \midrule
   \end{tabular*}
     \vspace*{-1.2em}
   }
 \end{table}


Table \ref{heur-avg} shows the overall performance impact of enabling heuristics over all instances that could be solved by at least one setting. On both sets, we see almost no change in total solving time. On \easy, the time to find the first solution decreases by $\reduction{0.1}{0.75}$ and the primal integral decreases by $\reduction{2037.41}{2351.8}$. The picture is slightly different on the numerically difficult test set. Here, the time to find the first solution decreases by $\reduction{1.30}{4.77}$, while the primal integral increases by $\percentage{422.9}{8670.7}$.

The worse performance and success rate on \numdiff is expected, considering that this test set was curated to contain instances with numerical challenges. On those instances floating-point heuristics find solutions that might either not be feasible in exact arithmetic or are not possible to fix for the repair heuristic.
In both test sets, the repair heuristic was able to find solutions, while not imposing any significant overhead in solving time.

\begin{table}[t]
  \caption{Comparison of new framework with and without primal heuristics (3~seeds, presolving enabled, instances where repair step was called)}
  \medskip
  \label{heur-avg}

  \small{
    \sisetup{round-mode=places}
   \begin{tabular*}{\textwidth}{l@{\extracolsep{\fill}}rS[table-format=3.1,round-precision=1,table-number-alignment=right]rS[table-format=5.1,round-precision=1,table-number-alignment=right]S[table-format=3.1,round-precision=1,table-number-alignment=right]rS[table-format=5.1,round-precision=1,table-number-alignment=right]}
   \toprule
   & & \multicolumn{3}{c}{heuristics disabled} & \multicolumn{3}{c}{heuristics enabled} \\
   \cmidrule(){3-5}\cmidrule(){6-8}
   {Test set} & {size} & {solv.\,time} & {time-to-first} & {primal int.} & {solv.\,time} & {time-to-first} & {primal int.} \\
   \midrule
   \easy & 82 & 32.50 &             0.75 &        2351.80 & 32.63 &             0.10 &        2037.19 \\
   \numdiff & 41 & 101.69 &             4.77 &        8670.74 & 103.06 &             1.30 &        9093.59 \\
   \midrule
   \end{tabular*}
     \vspace*{-1.2em}
   }
 \end{table}

\bigskip

\noindent
\textbf{Producing and verifying certificates.}
The possibility to log certificates as presented in~\cite{VIPR} is available in
the new framework and is extended to also work when the dual bounding method
bound-shift is active.
Presolving must currently be disabled, since \papilo does not yet support
generation of certificates.

Besides ensuring correctness of results, certificate generation is valuable to ensure correctness of the solver. Although it does not check the implementation itself, it can help identify and eliminate incorrect results that do not directly lead to fails. For example, on instance $x\_4$ from \numdiff, the original framework claimed infeasibility at the root node, and while the instance is indeed infeasible, we found the reasoning for this to be incorrect due to the use of a certificate.

Table \ref{vipr-cmp} reports the performance overhead when enabling certificates. Here we only consider instances that were solved to optimality by both versions since timeouts would bias the results in favor of the certificate.
We see an increase in solving time of $\increase{65.6}{32.6}$ on \easy and of $\increase{63.0}{41.6}$ on \numdiff. This confirms the measurements presented in \cite{VIPR}.
The increase is explained in part by the effort to keep track of the tree structure and print the exact dual multipliers, and in part by an increase in dual bounding time. The reason for the latter is that bound-shift by default only provides a safe objective value. The dual multipliers needed for the certificate must be computed in a postprocessing step, which introduces the overhead in safe bounding time. This overhead is larger on \easy, since bound-shift is called more often. 
The time spent in the verification of the certificate is on average significantly lower than the time spent in the solving process.
Overall, the overhead from printing and checking certificates is not negligible, but neither does it drastically impact the solvability of instances.

\begin{table}[t]
   \caption{Overhead for producing and verifying certificates on instances
     solved by both variants}
   \medskip

   \label{vipr-cmp}
   \small{
     \sisetup{round-mode=places,round-precision=1,table-number-alignment=right,table-format=2.1}
   \begin{tabular*}{\textwidth}{l@{\;\extracolsep{\fill}}rS[table-number-alignment=right]S[table-number-alignment=right]S[table-number-alignment=right]S[table-number-alignment=right]S[table-number-alignment=right]r}
   \toprule
   & & \multicolumn{2}{c}{certificate disabled} & \multicolumn{4}{c}{certificate enabled} \\
   \cmidrule(){3-4}\cmidrule(){5-8}
   {Test set} & {size} & {solving time} & {dbtime} & {solving time} & {dbtime} & {check time}  & {overhead} \\
   \midrule
   \easy & 53 & 32.64 &          9.13 &  65.62 &         16.54 & 0.88 &  103.9\,\% \\
   \numdiff & 21 & 41.55 &         11.86 & 63.00 &         17.97 &            0.45  & 52.6\,\%\\
   \midrule
   \end{tabular*}
     \vspace*{-1.2em}
   }
 \end{table}

\bigskip

\begin{table}[t]
   \caption{Comparison on \miplib~2017 benchmark set}
   \medskip

   \label{miplib-cmp}
   \small{
   \sisetup{round-mode=places,round-precision=1}
   \begin{tabular*}{\textwidth}{l@{\extracolsep{\fill}}r@{\extracolsep{\fill}}rrS[table-format=4.1]S[table-format=2.1]rrS[table-format=4.1]S[table-format=2.1]}
   \toprule
   && \multicolumn{4}{c}{original framework} & \multicolumn{4}{c}{new framework} \\
   \cmidrule(){3-6}\cmidrule(){7-10}
   {Test set} & {size} & {solved} & {found} & {time} & {gap} & {solved} & {found} & {time} & {gap} \\
   \midrule
   all & 240 & 17 & 74 & 6003.62 & {$\infty$} & 47 & 167 & 3928.05 & {$\infty$} \\ 
   both & 66 & 16 & 66 & 4180.96 & \SI{67.86}{\percent} & 29 & 66 & 1896.20 & \SI{33.8}{\percent} \\ 
   onesolved & 49 & 17 & 31 & 3317.47 & {$\infty$} & 47 & 47 & 505.13  & {$\infty$} \\ 
   \midrule
   \end{tabular*}
     \vspace*{-1.2em}
   }
 \end{table}


\noindent
\textbf{Performance comparison on MIBLIB~2017.}
As a final experiment, we wanted to evaluate the performance on a more ambitious and diverse test set.
To that end, we ran both the original framework and the revised framework with presolving and heuristics enabled on the recent \miplib~2017 benchmark set. The results in Table~\ref{miplib-cmp} show that the new framework solved $30$~instances more and the mean solving time decreased by $\reduction{505.1}{3317.5}$ on the subset ``onesolved'' of instances that could be solved to optimality by at least one solver. On more than twice as many instances at least one primal solution was found ($167$ vs. $74$). On the subset of $66$ instances that had a finite gap for both versions, the new algorithm achieved a gap of $33.8\%$ in arithmetic mean compared to $67.9\%$ in the original framework.

%
To conclude, we presented a substantially revised and extended solver for numerically exact
mixed integer optimization that significantly improves upon the existing state
of the art.
We also observe, however, that the performance gap to floating-point solvers is
still large.
This is not surprising, given that crucial techniques such as numerically safe
cutting plane separation, see, \eg \cite{CookDashSanjeebFukasawa2009}, are not yet included.
This must be addressed in
future research.

\bigskip

{\noindent\small\textbf{Acknowledgements.}
  We wish to thank Dan Steffy for valuable discussions on the revision of the original branch-and-bound framework, Leona Gottwald for creating \papilo, and Antonia Chmiela for help with implementing the primal repair heuristic.
}

\bibliographystyle{splncs04}
\bibliography{status-report}

\end{document}